\documentclass{article}%
\usepackage{amsmath}
\usepackage{amsfonts}
\usepackage{amssymb}
\usepackage{graphicx}%
\providecommand{\U}[1]{\protect\rule{.1in}{.1in}}

\begin{document}

\title{{\Large On Some Problems of Confidence Region Construction}}
\author{M. Evans, M. Liu, M. Moon, S. Sixta, S. Wei and S. Yang\\Dept. of Statistical\ Sciences\\University of Toronto}
\date{}
\maketitle

\begin{center}
\textbf{Abstract}
\end{center}

The general problem of constructing confidence regions is unsolved in the
sense that there is no algorithm that provides such a region with guaranteed
coverage for an arbitrary parameter $\psi\in\Psi.$ Moreover, even when such a
region exists, it may be absurd in the sense that either the set $\Psi$ or the
null set $\phi$ is reported with positive probability. An approach to the
construction of such regions with guaranteed coverage and which avoids
absurdity is applied here to several problems that have been discussed in the
recent literature and for which some standard approaches produce absurd
regions.\smallskip


\section{Introduction}

The confidence concept arises in statistics as follows: there is a statistical
model $\{f_{\theta}:\theta\in\Theta\}$ for data $x\in\mathcal{X},$ a marginal
parameter of interest $\psi=\Psi(\theta),$ where $\Psi:\Theta
\overset{onto}{\rightarrow}\Psi$ (with the same notation used here for the
function and its range), a desired confidence level $\gamma\in(0,1)$ and the
goal is to state a region $C(x)\subset\Psi$ such that $P_{\theta}(\Psi
(\theta)\in C(x))\geq\gamma$ for every $\theta\in\Theta.$ While there can be
different motivations for reporting such a region, the one considered here is
that there is an estimate $\psi(x)$ of the parameter of interest such that
$\psi(x)\in C(x)$ and the "size" of $C(x),$ together with the confidence
$\gamma,$ serve as an assessment of the accuracy of the recorded estimate. It
is well-known that confidence regions can sometimes give absurd answers as
discussed, for example, in Plante (2020). By absurd here is meant that $C(x)$
could be the null set or all of $\Psi$ with positive probability$,$ and so be
uninformative. In such situations it is difficult to see how reporting $C(x)$
can be regarded as a valid assessment of the accuracy of $\psi(x).$ Another
issue associated with confidence regions is that there isn't a theory that
prescribes how such a region can be constructed for a general problem.

The problem of error assessment via quoting a region $C(x)$, can also be
approached by adding a prior $\pi$ to the problem and providing a Bayesian
credible region having posterior content at least $\gamma.$ The Bayesian
approach has the virtues of the error assessment being based on the observed
data and such a region can always be constructed, say via the hpd (highest
posterior density) principle. There are criticisms that can be leveled at this
approach, however, as there is no assessment of the reliability of the
inference which is implicit in the frequentist approach via repeated sampling.
While the use of a prior is also sometimes criticized, the position taken here
is that this is no different than the use of a statistical model as, while the
model can be checked for its agreement with the observed data via model
checking, similarly a prior can be submitted to a check for prior-data
conflict, see Evans and Moshonov (2006), Evans (2015) and Nott et al. (2020).
There is also the issue of bias which is interpreted here as meaning that the
ingredients to the analysis, namely, the data collection procedure together
with the model and prior, can be chosen in such a fashion as to produce a
foregone conclusion. That such bias is possible is illustrated in Evans (2015)
and Evans and Guo (2021) where also a solution to this issue is developed.

Rather than invoke something like the hpd principle to construct a credible
region, the approach taken here is somewhat different. This is based on the
\textit{principle of evidence}: there is evidence in favor of a value $\psi$
if its posterior probability has increased over its prior probability,
evidence against $\psi$ if the posterior probability has decreased and there
is no evidence either way if they are equal. This simple principle has broad
implications not the least of which being that it makes little sense to allow
any reported region to include a value for which there is evidence against it
being true. In fact, a reported confidence or credible region can contain
$\psi$ values for which there is evidence against $\psi$ being the true value.
As such, it is more appropriate to quote what is called the \textit{plausible
region,} namely, those values of $\psi$ for which there is evidence in favor
of $\psi$ being true, see Evans (2015) and Section 2. The principle of
evidence also leads to a direct method for measuring and controlling bias
which comes in two forms for this problem. Here the \textit{implausible
region} refers to the set of $\psi$ values for which evidence against is
obtained.\smallskip

\noindent(i) \textit{Bias against} refers to the prior probability that the
plausible region does not contain the true value.

\noindent(ii) \textit{Bias in favor} refers to the prior probability that the
implausible region does not contain a meaningfully false value as defined in
Section 2.\smallskip

\noindent As discussed in Evans and Guo (2021), the control of bias is
equivalent to the a priori control of frequentist coverage probabilities for
the plausible region. Controlling bias against is equivalent to setting the
confidence of the plausible region, namely, the probability of it containing
the true value. Controlling bias in favor is typically equivalent to setting
the \textit{accuracy} of the plausible region where accuracy refers to the
probability of the plausible region covering false values. The measurement of
bias is reviewed in\ Section 2.

The end result of this approach is the best of both approaches to the problem,
namely, a Bayesian region with a particular posterior content that reflects
the uncertainty in the observed data, together with a guaranteed frequentist
confidence and accuracy, that reflects the reliability of the inference. The
reliability of an inference refers to the extent to which an inference is
trustworthy and, in general, Bayesian inferences do not address this issue. It
is important to note that these results hold for any proper prior and, at
least up to computational difficulties, can always be implemented. In
particular, there is no need to search for a prior that will provide an
appropriate confidence. So an elicited prior can be used, and moreover there
is no need for the posterior content and the confidence to agree, as they
refer to different aspects of the inference.

Section 2 discusses some necessary background and establishes the new result
that a plausible region is never absurd. Section 3 applies this approach to
several well-known problems where the construction of frequentist confidence
regions has proven to be at the very least difficult and, one could argue, for
which there is no current satisfactory solution. The methodology is general
and can be applied to any problem with a Bayesian formulation using proper
priors and so this provides a degree of unification between Bayes and frequentism.

\section{Relative Belief Inferences and Bias}

If the prior and posterior densities of $\Psi$ are denoted by $\pi_{\Psi}$ and
$\pi_{\Psi}(\cdot\,|\,x)$, then the relative belief ratio of $\Psi$ at $\psi
$\ is given by $RB_{\Psi}(\psi\,|\,x)=\pi_{\Psi}(\psi\,|\,x)/\pi_{\Psi}(\psi)$
and there is evidence in favor of $\psi$ when $RB_{\Psi}(\psi\,|\,x)>1,$
evidence against when $RB_{\Psi}(\psi\,|\,x)<1$ and no evidence either way
when $RB_{\Psi}(\psi\,|\,x)=1.$ This follows from the principle of evidence
when the prior distribution of $\Psi$ is discrete and follows via a limiting
argument in the general case, see Evans (2015). Actually, for much of the
discussion here, any valid measure of evidence can be used instead of the
relative belief ratio, where \textit{valid} means there is a cut-off that
determines evidence against versus evidence in favor according to the
principle of evidence. For example, a Bayes factor is a valid measure of
evidence also using the value 1 as the cut-off. As will be seen, the plausible
region and the measures of bias are independent of the valid measure of
evidence used so this is not an issue for the discussion here.

The set of values for which there is evidence in favor is the plausible region
$Pl_{\Psi}(x)=\{\psi:RB_{\Psi}(\psi\,|\,x)>1\}.$ When the $\psi$ values are
ordered by the amount of evidence via the relative belief ratio, the natural
estimate of $\Psi$ is given by $\psi(x)=\arg\sup_{\psi}RB_{\Psi}%
(\psi\,|\,x)\in Pl_{\Psi}(x).$ The posterior content of $Pl_{\Psi}(x)$
measures how strongly it is believed that the true value is in $Pl_{\Psi}(x)$
and the prior probability content of $Pl_{\Psi}(x)$ gives a measure of the
size or accuracy based on the observed data. So $\psi(x)$ can be considered a
highly accurate estimate when the posterior probability $\Pi_{\Psi}(Pl_{\Psi
}(x)\,|\,x)$ is high, as then there is a high degree of belief the true value
is in $Pl_{\Psi}(x),$ and the prior probability $\Pi_{\Psi}(Pl_{\Psi}(x))$ is
small, as then this set is small relative to the prior. Other measures of
posterior accuracy can also be quoted, such as the Euclidean measure or
cardinality measure of $Pl_{\Psi}(x)$ when relevant, but the posterior and
prior contents work universally for this purpose. Note that any other estimate
determined in this way from a valid measure of evidence will also lie in
$Pl_{\Psi}(x)$ and so produces no gain in accuracy over $\psi(x).$

It is possible, however, that there is bias in Bayesian inferences. For
example, suppose that the goal is to assess the hypothesis $H_{0}:\Psi
(\theta)=\psi_{0}.$ The relative belief ratio $RB_{\Psi}(\psi_{0}\,|\,x)$
indicates whether there is evidence in favor of or against $H_{0}$ and there
are several approaches to measuring the strength of this evidence but this is
not considered further here, see Evans (2015). Suppose that evidence against
$H_{0}$ is obtained but that there is a large prior probability of not getting
evidence in favor even when $H_{0}$ is true, namely, the probability
\begin{equation}
\text{bias against}_{\Psi}(\psi_{0})=M(RB_{\Psi}(\psi_{0}\,|\,x)\leq
1\,|\,\psi_{0}) \label{hyp1}%
\end{equation}
is large where $M(\cdot\,|\,\psi_{0})$ denotes the conditional prior
distribution of the data given that $H_{0}$ is true. It seems reasonable then
to treat the finding of evidence against $H_{0}$ as unreliable and it can be
said that there is an a priori \textit{bias against} $H_{0}.$ Similarly, using
a metric $d$\ on $\Psi,$ if evidence in favor of $H_{0}$ is obtained but
\begin{equation}
\text{bias in favor}_{\Psi}(\psi_{0})=\sup_{\psi:d(\psi_{0},\psi)\geq\delta
}M(RB_{\Psi}(\psi_{0}\,|\,x)\geq1\,|\,\psi) \label{hyp2}%
\end{equation}
is large, namely, there is a large prior probability of not obtaining evidence
against $H_{0}$ when it is \textit{meaningfully false}, as indicated by the
choice of the deviation $\delta,$ then it is said that there is \textit{bias
in favor} of $H_{0}.$ Note that $M(RB_{\Psi}(\psi_{0}\,|\,x)\geq1\,|\,\psi)$
generally decreases as $\psi$ moves away from $\psi_{0},$ so it is often only
necessary to consider values of $\psi$ satisfying $d(\psi_{0},\psi)=\delta$ to
determine the bias against. The value of $\delta$ is not arbitrary but is
determined by the application, as it represents the accuracy to which it is
desired to know the true value of $\psi$ which also determines the precision
of the measurement process that produces the data. Clearly there is some
similarity between the frequentist size and power of a test and the bias
against and bias in favor here but there is no suggestion that we are to
accept or reject $H_{0}.$ The purpose of the biases is to measure the
reliability of what the evidence in the observed data tells us about $H_{0}.$

The probability measures $M(\cdot\,|\,\,\psi)$ depend on the prior $\pi$ only
through the conditional prior $\pi(\cdot\,|\,\psi)\ $and do not depend on the
marginal prior $\pi_{\Psi}$ for the parameter of interest. As such the
probabilities determined by $M(\cdot\,|\,\,\psi)$ are essentially frequentist
in nature and similar to the use of distributions on parameters in mixed
models, namely, $\pi(\cdot\,|\,\psi)$ is used to integrate out nuisance
parameters. In fact, the bias probabilities (\ref{hyp1}) and (\ref{hyp2}) are
exactly frequentist but for the model given by $\{m(\cdot\,|\,\,\psi):\psi
\in\Psi\},$ where $m(\cdot\,|\,\,\psi)$ is the density of $M(\cdot
\,|\,\,\psi),$ and this corresponds to the original model when $\Psi
(\theta)=\theta.$

The average bias against a value of $\psi\sim\pi_{\Psi}$ can be written as
\begin{align}
\text{bias against}_{\Psi}  &  =E_{\Pi_{\Psi}}(M(RB_{\Psi}(\psi\,|\,x)\leq
1\,|\,\psi))=E_{\Pi_{\Psi}}(M(\psi\notin Pl_{\Psi}(x)\,|\,\psi))\nonumber\\
&  =1-E_{\Pi_{\Psi}}(M(\psi\in Pl_{\Psi}(x)\,|\,\psi)). \label{cov1}%
\end{align}
So (\ref{cov1}) is determined by the prior coverage probability $E_{\Pi_{\Psi
}}(M(\psi\in Pl_{\Psi}(x)\,|\,\psi))$ of the plausible region which will be
referred to hereafter as a (Bayesian) confidence as it is the prior
probability that $Pl_{\Psi}(x)$ contains the true value. Note that if an upper
bound can be obtained for $M(RB_{\Psi}(\psi\,|\,x)\leq1\,|\,\,\psi)$ as a
function of $\,\psi,$ then 1 minus this bound serves as a lower bound on the
confidence and, as will be seen, such a bound is commonly available. This
lower bound is then a confidence with respect to the model $\{m(\cdot
\,|\,\,\psi):\psi\in\Psi\}.$ Also the average bias in favor can be written as
\begin{align}
\text{bias in favor}_{\Psi}  &  =E_{\Pi_{\Psi}}\left(  \sup_{\psi_{\ast
}:d(\psi,\psi_{\ast})\geq\delta}M(RB_{\Psi}(\psi\,|\,x)\geq1\,|\,\psi_{\ast
})\right) \nonumber\\
&  =E_{\Pi_{\Psi}}\left(  \sup_{\psi_{\ast}:d(\psi,\psi_{\ast})\geq\delta
}M(\psi\notin\mbox{$Im_\Psi(x)$}\,|\,\psi_{\ast})\right)  , \label{cov2}%
\end{align}
which is the prior probability that a meaningfully false value is not in the
\textit{implausible region} $\mbox{$Im_\Psi(x)$}=\{\psi:RB_{\Psi}%
(\psi\,|\,x)<1\},$ the set of values for which there is evidence against. In
cases where the prior distribution of $\psi$ is continuous, then typically
(\ref{cov2}) is an upper bound on the prior probability of $Pl(x)$ covering a
meaningfully false value.

While it might be appealing to consider choosing the prior to make both these
biases small, this is the wrong approach as indeed experience indicates that
choosing a prior to minimize bias against simply increases bias in favor and
conversely. As discussed in Evans and Guo (2021), as the diffuseness of the
prior increases, typically bias in favor increases and bias against decreases.
The way to control these biases is, as established in Evans (2015), through
the amount of data collected as both biases converge to 0 as this increases.
As such, it is possible to control both the prior probability of $Pl_{\Psi
}(x)$ covering the true value and the prior probability of it covering a
meaningfully false value and so obtain a Bayesian inference with good
frequentist properties. Of course, this is similar to the use of coverage
probabilities in frequentist inference but the reported inferences are indeed
Bayesian while the biases are concerned with ensuring that the inferences are
reliable from a frequentist perspective.

A region $C$ for $\Psi$\ is called \textit{absurd} if it is possible that
$C(x)=\phi$ or $C(x)=\Psi$ with positive probability$.$ The following result
establishes that plausible regions can never be absurd in realistic
statistical contexts. The result can be viewed as a logical consistency result
for this approach to assessing the error in an estimate. For this let
$m(x)=\int_{\Theta}f_{\theta}(x)\,\Pi(d\theta)$ denote the prior predictive
density associated with the corresponding measure $M,m(x\,|\,\psi
)=\int_{\Theta}f_{\theta}(x)\,\Pi(d\theta\,|\,\psi)$ be the conditional prior
predictive density of the data given $\psi=\Psi(\theta)$ and put
\[
F=\{x:RB_{\Psi}(\psi\,|\,x)=1\text{ a.e. }\Pi_{\Psi}\}=\{x:m(x\,|\,\psi
)=m(x)\text{ a.e. }\Pi_{\Psi}\}
\]
where the last equality follows from the Savage-Dickey ratio result,\ namely,
$RB_{\Psi}(\psi\,|\,x)=m(x\,|\,\psi)/m(x).$ Note that, the conditional prior
distribution of the data given $F$ has no dependence on the parameter of
interest and, except in extraordinary circumstances, this set will have prior
probability 0, namely, $M(F)=0$. For, if $x\in F,$ then nothing can be learned
as there is no evidence in either direction for any value of $\psi.$\smallskip

\noindent\textbf{Theorem 1. }The plausible region for $\psi=\Psi(\theta)$ (i)
never satisfies $Pl_{\Psi}(x)=\Psi$ and (ii) satisfies $Pl_{\Psi}(x)=\phi$
with prior probability 0 when $M(F)=0.$

\noindent Proof: (i) Suppose that $Pl_{\Psi}(x)=\Psi.$ This is true iff
$RB_{\Psi}(\psi\,|\,x)>1$ for every $\psi\ $and so
\[
1<\int_{\Psi}RB_{\Psi}(\psi\,|\,x)\,\Pi_{\Psi}(d\psi)=\int_{\Psi}\frac
{\pi_{\Psi}(\psi\,|\,x)}{\pi_{\Psi}(\psi)}\,\Pi_{\Psi}(d\psi)=\int_{\Psi}%
\,\Pi_{\Psi}(d\psi\,|\,x)=1
\]
which is a contradiction. (ii) Now suppose $Pl_{\Psi}(x)=\phi,$ which is true
iff $RB_{\Psi}(\psi\,|\,x)\leq1$ for every $\psi.$ Since $M(F)=0,$ this
implies that for any $x\notin F,$ the set $A(x)=\{\psi:RB_{\Psi}%
(\psi\,|\,x)=1\}$ has $\Pi_{\Psi}(A(x))<1$ which implies $0<1-\Pi_{\Psi
}(A(x))=\Pi_{\Psi}(\mbox{$Im_\Psi(x)$}).$ Then
\begin{align*}
1  &  =\Pi_{\Psi}(A(x))+\Pi_{\Psi}(\mbox{$Im_\Psi(x)$})\\
&  >\int_{A(x)}RB_{\Psi}(\psi\,|\,x)\,\Pi_{\Psi}(d\psi)+\int%
_{\mbox{$Im_\Psi(x)$}}RB_{\Psi}(\psi\,|\,x)\,\Pi_{\Psi}(d\psi)\\
&  =\int_{\Psi}RB_{\Psi}(\psi\,|\,x)\,\Pi_{\Psi}(d\psi)=\int_{\Psi}\,\Pi
_{\Psi}(d\psi\,|\,x)=1
\end{align*}
which is a contradiction. \smallskip

It is also possible to construct credible regions based on the relative belief
ratio as in $C_{\gamma}(x)=\{\psi:RB_{\Psi}(\psi\,|\,x)\geq r_{\gamma}\}$
where $r_{\gamma}=\sup\{r:\Pi_{\Psi}(RB_{\Psi}(\psi\,|\,x)<r\,|\,x)\leq
1-\gamma\}$ as then $\Pi_{\Psi}(C_{\gamma}(x)\,|\,x)\geq\gamma.$ As with all
relative belief inferences, the relative belief credible regions are invariant
under smooth reparameterizations while hpd regions are not. This means that
the computation of a $\gamma$-relative belief region can be carried out in any
parameterization while each parameterization leads to a potentially different
hpd credible region. With both approaches, however, it is impossible to say a
priori that all the elements of the region will have evidence in their favor.
For relative belief regions, however, it is guaranteed that for any\ $\gamma
\leq\Pi_{\Psi}(Pl_{\Psi}(x)\,|\,x)$ then $C_{\gamma}(x)\subset Pl_{\Psi}(x)$
and there is evidence in favor of each element of $C_{\gamma}(x)$ so such a
region can be also be reported. There are also a variety of optimality
properties satisfied by relative belief credible regions, see Evans (2015).
The property of importance for the discussion here, however, is that for the
plausible region it can be determined a priori how much data to collect to
ensure appropriate coverage probabilities and that doesn't seem to be
available for a credible region in general.

It is also the case, as established in Evans and Guo (2021), that plausible
regions possess additional good, and even optimal, properties beyond those
already cited like parameterization invariance and no dependence on the valid
measure of evidence used. For example, the prior probability of $Pl_{\Psi}(x)$
covering the true value is always greater than or equal to the prior
probability of $Pl_{\Psi}(x)$ covering a false value which in frequentist
theory is known as the unbiasedness property for confidence regions. As an
example of an optimal property, when the prior $\Pi_{\Psi}$ is continuous,
then among all regions $C$ satisfying $M(\psi\in C(X)\,|\,\psi)\geq M(\psi\in
Pl_{\Psi}(X)\,|\,\psi)$ for every $\psi,$ namely, the conditional prior
probability that $C$ covers the true value is as large as this probability for
$Pl_{\Psi},$ then $Pl_{\Psi}$ maximizes the prior probability of not covering
a false value and there is a similar optimality property for the discrete
case. The implication of this is that, if one considers another way of
expressing evidence that leads to the region $C,$ then provided its coverage
probabilities are as large as those of $Pl_{\Psi},$ as otherwise it presumably
wouldn't be considered, then $C$ cannot do better than $Pl_{\Psi}$ with
respect to accuracy. This is really an optimality property for the principle
of evidence and there are other such results.

\section{Examples}

There are a variety of problems discussed in the literature where issues
concerning either absurd confidence regions are obtained or it is unclear how
to construct a $\gamma$-confidence region for a general parameter $\psi
=\Psi(\theta)$. The following examples show that the approach via the
principle of evidence can deal successfully with such problems.

\subsection{Fieller's Problem}

This is a well-known problem, as discussed in Geary (1930), Fieller (1954),
Hinkley (1969) and more recently in Pham-Gia et al. (2006) and Ghosh et al.
(2006) where a wide range of applications are noted. Ghosh et al. (2006) is
concerned with confidence intervals for ratios of regression coefficients in a
normal linear model and it is shown that certain integrated likelihoods do not
produce absurd intervals and this is now a consequence of the general Theorem
1. This problem is also discussed in Fraser et al. (2018) where it appears as
problems A and B of a set of problems for inference proposed by D. R. Cox.

For this there are two samples $x=(x_{1},\ldots,x_{m})$ $i.i.d$. $N(\mu
,\sigma_{0}^{2})$ independent of $y=(y_{1},\ldots,y_{n})$ $i.i.d.$
$N(\nu,\sigma_{0}^{2})$ where $(\mu,\nu)\in%
\mathbb{R}
^{2}$ is unknown. So it is supposed that the means are unknown but the
variances are known and common. The discussion can be generalized to allow for
unknown variances as well, with no changes to the basic results, but the
essential problem arises in the simpler context. The problem then is to make
inference about the ratio of means $\psi=\Psi(\mu,\nu)=\mu/\nu$ and, in
particular, construct a confidence interval for this quantity. It is assumed
here that model checking has not led to any suspicions concerning the validity
of the models. As such the data can be reduced to the minimal sufficient
statistic $(\bar{x},\bar{y})$ where $\bar{x}\sim N(\mu,\sigma_{0}^{2}/n)$
independent of $\bar{y}\sim N(\nu,\sigma_{0}^{2}/m).$

Confidence regions for $\psi$ can be obtained via a pivotal statistic given by
$(\bar{x}-\bar{y}\psi)/\sigma_{0}\sqrt{1/m+\psi^{2}/n}\sim N(0,1)$ but this
can produce absurd regions.\ For example, if a $\gamma$-confidence interval is
required for $\psi$ then, with $z_{p}$ denoting the $p$-th quantile of a
$N(0,1),$ the region equals $%
\mathbb{R}
^{1}$ whenever $|\bar{y}|<\sigma_{0}z_{(1+\gamma)/2}/n^{1/2}$ and $m\bar
{x}^{2}+n\bar{y}^{2}<z_{(1+\gamma)/2}^{2}\sigma_{0}^{2}.$ Sometimes the region
can be a so-called \textit{exclusive region} of the form $(-\infty,a(\bar
{x},\bar{y}))\cup(b(\bar{x},\bar{y}),\infty)$ with $a(\bar{x},\bar{y}%
)<b(\bar{x},\bar{y}).$ While an interval might be preferred, there is nothing
illogical about an exclusive region as can be seen by considering the $\gamma
$-confidence interval $\bar{x}\pm\sigma_{0}z_{(1+\gamma)/2}/m^{1/2}$ for
$\mu.$ If this interval includes $0,$ then necessarily the $\gamma$-confidence
region for $1/\mu$ has the exclusive form $(-\infty,1/(\bar{x}-\sigma
_{0}z_{(1+\gamma)/2}))\cup(1/(\bar{x}+\sigma_{0}z_{(1+\gamma)/2}),\infty).$
The same reasoning applies in Fieller's problem and one can always
reparameterize by making inference instead about $\psi^{-1}$ to obtain an
interval. The problem of exclusive regions is a consequence of the
parameterization but that is not the case with absurd regions as this
represents a defect in the inference.

The relative belief approach requires the specification of a prior and for
this conjugate priors $\mu\sim N(\mu_{0},\tau_{10}^{2})$ independent of
$\nu\sim N(\nu_{0},\tau_{20}^{2}),$ will be used. This requires an elicitation
for the quantities $(\mu_{0},\tau_{10}^{2},\nu_{0},\tau_{20}^{2})$ which can
proceed as follows. First specify $(m_{1},m_{2})$ such that the true value of
$\mu\in(m_{1},m_{2})$ with virtual certainty, say with prior probability
$\gamma=0.99.$\textbf{ }Then put $\mu_{0}=(m_{1}+m_{2})/2$ and solve
$\Phi((m_{2}-\mu_{0})/\tau_{10})-\Phi((m_{1}-\mu_{0})/\tau_{10})=\gamma$ for
$\tau_{10},$ so the prior on $\mu$ is now determined. This step could also be
applied to obtain the prior\ for $\nu$ but it is supposed instead that there
is information about the true value of $\psi$ expressed as $\psi\in
(r_{1},r_{2})$ with virtual certainty for fixed constants $r_{1}<r_{2}$. A
value $\psi_{0}\in(r_{1},r_{2})$ is then selected, which could be a
hypothesized value for this quantity or just the central value, and then take
$\nu_{0}=\mu_{0}/\psi_{0}.$ Finally, requiring $\nu\in(m_{1}/r_{2},m_{2}%
/r_{1})$ with virtual certainty determines $\tau_{20}$ via $\Phi((m_{2}%
/r_{1}-\nu_{0})/\tau_{20})-\Phi((m_{1}/r_{2}-\nu_{0})/\tau_{20})=\gamma$ and
this gives the prior for $\nu$\textbf{.} This is just one method for eliciting
the prior and an alternative could be more suitable in a given application.
Once a prior has been determined and the data obtained, the prior is subjected
to a check for prior-data conflict and it is assumed here that the prior has
passed such a check.

Some numerical examples are carried along for illustration
purposes.\textbf{\smallskip}

\noindent\textbf{Example 1. }\textit{Simulation example (the data, model and
prior).}

Suppose $m=n=10,\mu=20,\nu=10,\sigma_{0}^{2}=1$ so the true value is
$\psi=20/10=2\mathbf{.}$\textbf{ }Data was generated leading to the mss
$(\bar{x},\bar{y})=(20.188,10.699).$ For the prior elicitation suppose
$(m_{1},m_{2})=(10,25),$ so $\mu_{0}=17.5,\tau_{10}=2.912$ and with
$(r_{1},r_{2})=(1,3),\psi_{0}=2\,$\ then $\nu_{0}=8.75,\tau_{20}=2.336.$ The
value $\psi_{0}=2$ is chosen as the hypothesis $H_{0}:\Psi(\mu,\nu)=2$ will be
subsequently assessed to see how the approach performs with a true hypothesis.
Inverting the pivotal leads to the $0.95$-confidence region $(1.770,2.016)$
for $\psi$ which just includes the true value.\textbf{\smallskip}

\noindent\textbf{Example 2. }\textit{Cox's examples (the data, model and
prior).}

For the Cox A problem, $m=n,\sigma_{0}^{2}/n=1,(\bar{x},\bar{y})=(10,0.5)$
which produces the exclusive $0.95$-confidence region $(-\infty,-6.752)\cup
(3.968,\infty)$ via the pivotal. For Cox B the only change is that now
$(\bar{x},\bar{y})=(0.5,0.5)$ and the $0.95$-confidence region is $%
\mathbb{R}
^{1}$ and so is absurd. No priors were prescribed for either problem, so here
we take fairly noninformative priors that avoid prior-data conflict. For
problem A suppose $\mu\sim N(12,3)$ independent of $\nu\sim N(0,3)$ and for
problem B suppose both priors are $N(0,3).$ $\blacksquare$\textbf{\smallskip}

Putting%
\begin{equation}
\tau_{20}^{2}(\psi)=\left(  \psi^{2}/\tau_{10}^{2}+1/\tau_{20}^{2}\right)
^{-1},\text{ }\nu_{0}(\psi)=\tau_{20}^{2}(\psi)\left(  \psi\mu_{0}/\tau
_{10}^{2}+\nu_{0}/\tau_{20}^{2}\right)  , \label{form1}%
\end{equation}
then the exact prior density of $\psi$ is
\begin{align*}
\pi_{\Psi}(\psi)  &  =\frac{2\tau_{20}^{2}(\psi)}{\sqrt{2\pi}\tau_{10}%
\tau_{20}}\exp\left\{  -\frac{1}{2}\frac{\tau_{20}^{2}(\psi)(\mu_{0}-\nu
_{0}\psi)^{2}}{\tau_{10}^{2}\tau_{20}^{2}}\right\}  \times\\
&  \left\{  \varphi\left(  \frac{\nu_{0}(\psi)}{\tau_{20}(\psi)}\right)
+\frac{\nu_{0}(\psi)}{\tau_{20}(\psi)}\Phi\left(  \frac{\nu_{0}(\psi)}%
{\tau_{20}(\psi)}\right)  -\frac{1}{2}\right\}  .
\end{align*}
Note that when $\mu_{0}=\nu_{0}=0$ then $\pi_{\Psi}$ is a (rescaled) Cauchy
density so in general this distribution has quite long tails. The same formula
works for the posterior with substitutions as in (\ref{form1}) since
\begin{align*}
\mu\,|\,\bar{x}  &  \sim N(\left(  \mu(\bar{x}),\left(  \frac{m}{\sigma
_{0}^{2}}+\frac{1}{\tau_{10}^{2}}\right)  ^{-1}\right)  \text{ with }\mu
(\bar{x})=\left(  \frac{m}{\sigma_{0}^{2}}+\frac{1}{\tau_{10}^{2}}\right)
^{-1}\left(  \frac{m\bar{x}}{\sigma_{0}^{2}}+\frac{\mu_{0}}{\tau_{10}^{2}%
}\right) \\
&  \text{independent of}\\
\nu\,|\,\bar{y}  &  \sim N\left(  \nu(\bar{y}),\left(  \frac{n}{\sigma_{0}%
^{2}}+\frac{1}{\tau_{20}^{2}}\right)  ^{-1}\right)  \text{ with }\nu(\bar
{y})=\left(  \frac{n}{\sigma_{0}^{2}}+\frac{1}{\tau_{20}^{2}}\right)
^{-1}\left(  \frac{n\bar{y}}{\sigma_{0}^{2}}+\frac{\nu_{0}}{\tau_{20}^{2}%
}\right)  .
\end{align*}
So the relative belief ratio $RB_{\Psi}(\psi\,|\,\bar{x},\bar{y})$ is
available in closed form.

For a general problem, a closed form is typically not available for the prior
and posterior densities of marginal parameters of interest. In an application,
however, there is a difference $\delta>0$ that represents the accuracy with
which it is desired to know the true value. This quantity is a major input
into sample size considerations. The approach then is to partition the
effective prior range of $\psi,$ as determined via a simulation from the prior
of $(\mu,\nu),$ into subintervals of length $\delta$ with the midpoint of each
interval taken as representative of the values in that subinterval. The prior
and posterior contents of these subintervals are determined via a simulation
and then density histograms are used to approximate $\pi_{\Psi}(\psi)$ and
$\pi_{\Psi}(\psi\,|\,\bar{x},\bar{y})$ which in turn gives an approximation to
$RB_{\Psi}(\psi\,|\,\bar{x},\bar{y})$ that can be used to determine the
inferences.\textbf{\smallskip}

\noindent\textbf{Example 1. }\textit{Simulation example (the inferences).}

\begin{figure}[tb]%
\centering
\includegraphics[
height=3.1912in,
width=3.1998in
]%
{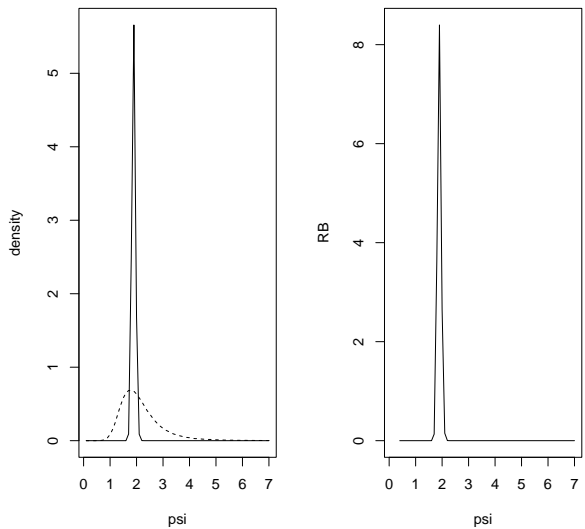}%
\caption{Plots (left panel) of the prior (- - -) and the posterior (--)
densities and (right panel) of the relative belief ratio of $\psi$ in Example
1.}%
\label{fig1}%
\end{figure}

The above approximation procedure was carried out, using the values recorded
when $\delta=0.1$ was chosen for the accuracy. Figure \ref{fig1} provides
plots of $\pi_{\Psi}$ and $\pi_{\Psi}(\cdot\,|\,\bar{x},\bar{y})$ and
$RB_{\Psi}(\cdot\,|\,\bar{x},\bar{y}).$ Due to the long-tailed feature of the
prior some extreme values of $\psi$ are obtained and this is reflected in the
range over which these distribution have been plotted. Relatively smooth
estimates are obtained based on Monte Carlo sample sizes of $N=10^{5}$ and
these can be seen to closely approximate the true functions. One approach for
coping with the long-tail is to calculate the ecdf $\hat{F}_{\Psi}$ of $\Psi$
based on a large simulation sample and take $(\psi_{\min},\psi_{\max}%
)=(\hat{F}_{\Psi}^{-1}(0.0005),\hat{F}_{\Psi}^{-1}(0.9995))$ so this ignores
$0.001$ of the probability in the tails which is what was done here. Another
possibility, which avoids the truncation, is to transform to $\omega=G(\psi)$
where $G$ is a long-tailed cdf like a Cauchy (or even sub-Cauchy) and
transform the initial partition to $(G(\psi_{\min}-\delta/2),G(\psi_{\min
}+\delta/2)],\ldots,(G(\psi_{\max}-\delta/2),G(\psi_{\max}+\delta/2)].$ All
inferences for $\psi$ can then be obtained from those for $\omega$\ via the
transformation $\psi=G^{-1}(\omega)$ due to the invariance of relative belief
inferences under reparameterizations.

The relative belief estimate is given by $\psi(x)=1.90$ with plausible region
$Pl_{\Psi}(x,y)=(1.75,2.05)$ having posterior content $0.982$ and prior
content $0.200.$ So the plausible region contains the true value, and note
that the estimate is reasonably accurate for a relatively small amount of
data. $\blacksquare$\textbf{\smallskip}%

\noindent\textbf{Example 2. }\textit{Cox's examples (the inferences).}

For the Cox A problem, the plausible region is $(-\infty,-17.60)\cup
(6.40,\infty)$ having posterior content $0.791$ and prior content $0.515.$ So
the inferences are not very precise. For the Cox B problem, the plausible
region is $(0.10,\infty)$ having posterior content $0.534$ and prior content
$0.498$ and the absurd interval is avoided. Both cases can be considered
extreme as there is little data relative to the variance $\sigma_{0}^{2}%
.$\textbf{\smallskip}

Now consider the bias calculations. To compute the biases for hypothesis
assessment it is necessary to compute%
\begin{equation}
M(RB_{\Psi}(\psi_{0}\,|\,\bar{x},\bar{y})\leq1\,|\,\psi) \label{bias1}%
\end{equation}
for various values\ of $\psi$ where $(\bar{x},\bar{y})$\ is generated from the
conditional prior predictive given $\psi$ and to compute the biases for
estimation we need to be able to compute (\ref{bias1}) for values of $\psi
_{0}\sim\pi_{\Psi}$ and then average.\textbf{ }So it is necessary to: (i)
generate $(\bar{x},\bar{y})$\ from its conditional prior predictive
$M(\cdot\,|\,\psi)$ and (ii) compute $RB_{\Psi}(\psi\,|\,\bar{x},\bar{y})$ and
compare it to 1.

For (i) the following sequential algorithm will work:\smallskip

\noindent1. generate $\nu\,|\,\psi\,\sim\pi(\cdot\,|\,\psi),$ 2. generate
$\bar{y}\,|\,(\psi,\nu)\sim N(\nu,\sigma_{0}^{2}/n),$ 3. generate $\bar
{x}\,|\,(\psi,\nu,\bar{y})\sim N(\psi\nu,\sigma_{0}^{2}/m).$\smallskip

\noindent Steps 2 and 3 are straightforward while step 1 requires the
development of a suitable algorithm. The joint prior density of $(\psi
,\nu,\bar{x},\bar{y})$ is proportional to
\[
|\nu|\exp\left\{  -\frac{1}{2}\left[  \frac{m(\bar{x}-\psi\nu)^{2}+n(\bar
{y}-\nu)^{2}}{\sigma_{0}^{2}}+\frac{(\psi\nu-\mu_{0})^{2}}{\tau_{10}^{2}%
}+\frac{(\nu-\nu_{0})^{2}}{\tau_{20}^{2}}\right]  \right\}
\]
\noindent which implies that $\pi(\nu\,|\,\psi)\propto|\nu|\exp\left\{
-\left(  \nu-\nu_{0}(\psi)\right)  ^{2}/2\tau_{20}^{2}(\psi)\right\}  $
where\ $\tau_{20}^{2}(\psi)$ and $\nu_{0}(\psi)$ are as specified in
(\ref{form1}). Therefore, $\pi(\cdot\,|\,\psi)$ is close to a normal density
but for the factor $|\nu|.$ Transforming $\nu\rightarrow z=\left(  \nu-\nu
_{0}(\psi)\right)  /\tau_{20}(\psi)$ we need to be able to generate $z$ from a
density of the form $g(z)\propto|\nu_{0}(\psi)+\tau_{20}(\psi)z|\varphi
(z)\propto|z-z_{0}|\varphi(z)$ where $z_{0}=-\nu_{0}(\psi)/\tau_{20}(\psi
).$\textbf{ }Using $d\varphi(z)/dz=-z\varphi(z),$ then $\int_{-\infty}^{z_{0}%
}-(z-z_{0})\varphi(z)\,dz=z_{0}\Phi(z_{0})+\varphi(z_{0})$ and $\int_{z_{0}%
}^{\infty}(z-z_{0})\varphi(z)\,dz=-z_{0}(1-\Phi(z_{0}))+\varphi(z_{0})$ which
implies
\begin{align*}
g(z)  &  =p(z_{0})I_{(-\infty,z_{0}]}(z)g_{1}(z)+(1-p(z_{0}))I_{(z_{0}%
,\infty)}(z)g_{0}(z),\text{ where}\\
g_{1}(z)  &  =\frac{(z_{0}-z)\varphi(z)}{z_{0}\Phi(z_{0})+\varphi(z_{0}%
)}\text{ when }z\leq z_{0},\\
g_{0}(z)  &  =\frac{(z-z_{0})\varphi(z)}{-z_{0}(1-\Phi(z_{0}))+\varphi(z_{0}%
)}\text{ when }z>z_{0},\text{ and }\\
p(z_{0})  &  =\frac{z_{0}\Phi(z_{0})+\varphi(z_{0})}{-z_{0}+2(z_{0}\Phi
(z_{0})+\varphi(z_{0}))}.
\end{align*}
So with probability $p(z_{0}),$ generate $z$ from $g_{1}$ and otherwise
generate $z$ from $g_{0}\mathbf{.}$\textbf{ }The cdf of $g_{1}$ for $z\leq
z_{0}$ equals%
\[
G_{1}(z)=\int_{-\infty}^{z}g_{1}(x)\,dx=\frac{z_{0}\Phi(z)+\varphi(z)}%
{z_{0}\Phi(z_{0})+\varphi(z_{0})}%
\]
and to generate from $g_{1}$ via inversion generate $u\sim U(0,1)$ and solve
$G_{1}(z)=u$ for $z$ by bisection. To start the bisection set $z_{up}=z_{0}$
and for some $\epsilon>0$ iteratively evaluate $G_{1}(-i|z_{0}-\epsilon|)$ for
$i=0,1,\ldots$ until $G_{1}(-i|z_{0}-\epsilon|)\leq u$ setting $z_{low}%
=-i|z_{0}-\epsilon|$ as this guarantees $G_{1}(z_{low})\leq u\leq G_{1}%
(z_{up})=1$ so bisection will work.\textbf{ }The cdf of $g_{0}$ is for
$z>z_{0}$
\[
G_{0}(z)=\int_{z_{0}}^{z}g_{0}(x)\,dx=\frac{z_{0}(\Phi(z)-\Phi(z_{0}%
))+(\varphi(z)-\varphi(z_{0}))}{z_{0}(1-\Phi(z_{0}))-\varphi(z_{0})}%
\]
and for this start bisection with $z_{low}=z_{0}$ and iteratively evaluate
$G_{0}(i|z_{0}+\epsilon|)$ until $u\leq G_{0}(i|z_{0}+\epsilon|)$ setting
$z_{up}=i|z_{0}+\epsilon|$ so bisection will work. Finally, when $z$ is
obtained put $\nu=\nu(\psi)+\tau_{20}(\psi)z$ to get the appropriately
generated value of $\nu.$\textbf{ }An interesting consequence of this
algorithm is that it must be true that $0<-z(1-\Phi(z))+\varphi(z)$ for every
$z$ and this implies the well-known \textit{Mills ratio inequality }%
$(1-\Phi(z))/\varphi(z)<1/z$ when $z\geq0$ and $\Phi(z)/\varphi(z)<1/|z|$ when
$z\leq0$\textit{ }which gives useful bounds on tail probabilities for the
normal distribution when $|z|$ is large.

To determine (\ref{bias1}) the value $RB_{\Psi}(\psi_{0}\,|\,\bar{x},\bar{y})$
needs to be computed for each\ generated value of $(\bar{x},\bar{y}).$ This
can be carried out as previously using the discretized version but using the
closed form version is much more efficient. It might seem more appropriate to
use the exact form also for inferences but, because we wish to incorporate the
meaningful difference $\delta$ into the inferences, the discretized version is
much more efficient for those computations. Note too that a high degree of
accuracy is not required for the bias computations.

Now consider the biases in the numerical problems being
considered.\textbf{\smallskip}

\noindent\textbf{Example 1. }\textit{Simulation example (the biases).}

\begin{figure}[tb]%
\centering
\includegraphics[
height=2.5079in,
width=2.514in
]%
{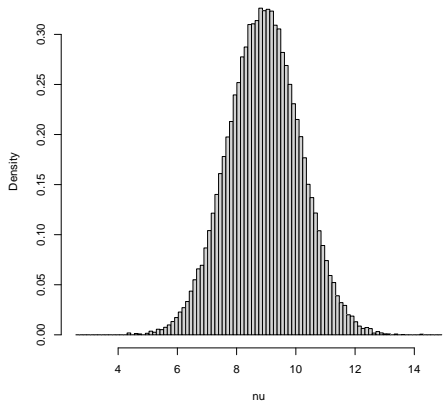}%
\caption{Density histogram of a sample of $10^{5}$ from the conditional prior
of $\nu$ given $\psi_{0}=2.$}%
\label{nusamp}%
\end{figure}

The hypothesis assessment problem is $H_{0}:\psi_{0}=2$ then, using the
elicited values of $(\mu_{0},\tau_{10}^{2},\nu_{0},\tau_{20}^{2}),$ leads to
$\tau_{20}^{2}(2)=1.5239,\nu_{0}(2)=8.7502$ and $z_{0}=$ $-7.0883.$ Figure
\ref{nusamp} is density histogram of a sample of $10^{5}$ from $\pi
(\cdot\,|\,\psi_{0}).$%

\begin{figure}[tb]%
\centering
\includegraphics[
height=2.4275in,
width=2.4344in
]%
{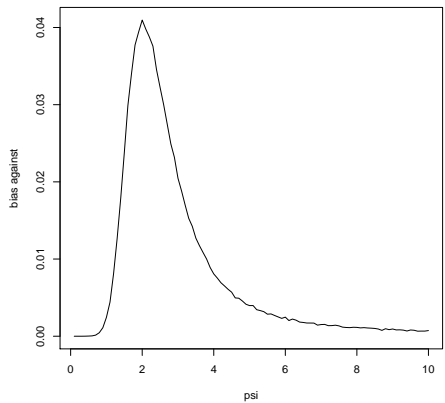}%
\caption{The bias against as a function of $\psi$ in Example 1.}%
\label{biasagex1}%
\end{figure}

To get the bias against $H_{0},$ use the sequential algorithm to generate
$(\bar{x},\bar{y})\sim M(\cdot\,|\,\psi_{0})$, compute and compare $RB_{\Psi
}(\psi_{0}\,|\,\bar{x},\bar{y})$ to $1$ for a large number of repetitions
recording the proportion of times $RB_{\Psi}(\psi_{0}\,|\,\bar{x},\bar{y}%
)\leq1.$ In this problem the value $M(RB_{\Psi}(\psi_{0}\,|\,\bar{x},\bar
{y})\leq1\,|\,\psi_{0})=0.04$ was obtained based on a Monte Carlo sample of
$10^{5}$ and so there is no real bias against $H_{0}:\psi_{0}=2.$ Figure
\ref{biasagex1} is a plot of $M(RB_{\Psi}(\psi\,|\,\bar{x},\bar{y}%
)\leq1\,|\,\psi)$ versus $\psi$ which is maximized at $\psi=2.0$ and takes the
value $0.040$ there. This implies that the conditional prior probability the
plausible plausible region contains the true value is at least $0.960$ for all
$\psi$\ and so can be considered as a $0.96$-confidence interval for $\psi.$
If instead we had $m=n=20,$ then the maximum bias against is $0.028$ and the
plausible region would then be $0.972$-confidence interval for $\psi$ and of
course larger sample sizes will just increase the confidence.

\begin{figure}[tb]%
\centering
\includegraphics[
height=2.3065in,
width=2.3125in
]%
{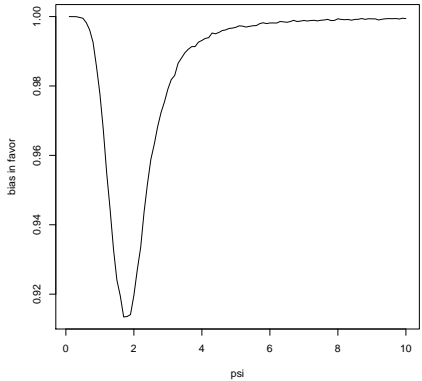}%
\caption{Bias in favor as a function of $\psi$ in Example 1.}%
\label{biasinfavex1}%
\end{figure}

To get the bias in favor of $H_{0},$ use the sequential algorithm to generate
$(\bar{x},\bar{y})\sim M(\cdot\,|\,\psi_{0}+\delta/2)$, compute and compare
$RB_{\Psi}(\psi_{0}\,|\,\bar{x},\bar{y})$ to 1, for a large number of
repetitions record the proportion of times $RB_{\Psi}(\psi_{0}\,|\,\bar
{x},\bar{y})\geq1,$ and also do this for $(\bar{x},\bar{y})\sim M(\cdot
\,|\,\psi_{0}-\delta/2)$ and the maximum of the two is an upper bound on the
bias in favor. In this case the value $0.92\,$is obtained which is very high
indicating that there is substantial bias in favor of the hypothesis. In other
words, there is a substantial prior probability that evidence in favor of the
hypothesis will be obtained even when it is meaningfully false as determined
by $\delta.$ Of course, sample size is playing a role here as well as
$\delta.$ For $m=n=20$ the upper bound equals $0.91$, for $m=n=100$ the upper
bound equals $0.71$, while for $m=n=500$ the upper bound equals $0.09.$ So
$m=n=10$ is not enough data to ensure that evidence in favor of $H_{0}$ will
not be obtained when it is meaningfully false with $\delta=0.1$ and more data
needs to be collected to avoid this. For the bias in favor for estimation a
sample of $\psi\sim\pi_{\Psi}$ values is generated and the bias in favor of
$\psi$ at $\psi\pm\delta/2$ is determined and then averaged. Figure
\ref{biasinfavex1} is a plot of the bias in favor as a function of $\psi$ and
the average value is $0.94$ which is an upper bound on the the prior
probability that the plausible region contains a meaningfully false value.
When $m=n=20$ the upper bound equals $0.92,$ when $m=n=100$ the upper bound
equals $0.69$ and when $m=n=500$ the upper bound equals $0.26.$ The value of
$\delta$ is determined by the application and taking it too small clearly
results in the requirement of overly large sample sizes to get the bias in
favor small. For example, with $m=n=10$ and $\delta=0.5,$ then the bias in
favor for estimation is $0.33$ while for $\delta=1.0$ it is $0.12,$ and with
$m=n=20$ these values are $0.21$ and $0.07$, respectively.\textbf{\smallskip}%

\noindent\textbf{Example 2. }\textit{Cox's examples (the biases).}

For the first problem an upper bound on the bias against is given by $0.18$ so
the coverage probability for the plausible region is at least $0.82.$ For the
second problem an upper bound on the bias against is given by $0.24$ so the
coverage probability for the plausible region is at least $0.76.$ These
coverages are quite reasonable given the small sample sizes relative to
$\sigma_{0}^{2}.$

\subsection{Mandelkern's Examples}

Mandelkern (2002) discusses several problems in physics where confidence
intervals are required but for which no acceptable solution exists. These are
problems where standard statistical models are used and, in the unconstrained
case, well-known confidence intervals are available for $\psi$ but physical
theory demands that the true value lie in $\Psi_{0}\subset\Psi$ for a proper
subset $\Psi_{0}.$ If $C(x)$ is a $\gamma$-confidence region for unconstrained
$\psi,$ then it is certainly the case that $C(x)\cap\Psi_{0}$ is a $\gamma
$-confidence region under the constraint. While this has the correct coverage
probability, however, in general $C(x)\cap\Psi_{0}$ can equal $\phi$ with
positive probability and so this solution is absurd. As is now demonstrated
the approach discussed here provides an effective solution to this problem.

Mandelkern's examples are now described together with the
solutions.\textbf{\smallskip}

\noindent\textbf{Example 3.} \textit{Location-normal with constrained mean.}

The model here is that a sample $x=(x_{1},\ldots,x_{n})$ has been obtained
from a distribution in $\{N(\mu,\sigma_{0}^{2}):\mu\in(l_{0},u_{0})\}$ where
$\sigma_{0}^{2}$ is known and $\mu$ is known to lie in the interval
$(l_{0},u_{0})$ where $l_{0},u_{0}\in%
\mathbb{R}
^{1}\cup\{-\infty,\infty\}$ with $l_{0}<u_{0}.$ Mandelkern discusses
inferences concerning the mass $\mu$ of a neutrino so \ in that case
$l_{0}=0.$ The measurements are taken to a certain accuracy and this is
reflected in the specification of the quantity $\delta$ which is the accuracy
to which it is desired to know $\mu$ which may indeed be larger than the
accuracy of the measurements. This leads to a grid of possible values for
$\mu,$ say $\mu_{1}<\mu_{2}<\cdots$ and such that $\left\vert \mu_{i}%
-\mu_{i+1}\right\vert =\delta.$ So when $l_{0}$ and $u_{0}$ are both finite
the possible values of $\mu$ are given by $\mu_{i}=l_{0}+(i-1/2)\delta$ for
$i=1,\ldots,(u_{0}-l_{0})/\delta,$ and it is supposed that these values are
such that $(u_{0}-l_{0})/\delta$ is an integer. It is certainly possible for
one or both of $l_{0},u_{0}$ to be infinite but typically there are lower and
upper bounds on what a measurement can equal. So in practice a finite number
of such intervals with possibly two tail intervals, which contain very little
prior probability, suffices. For example, consider measuring a length to the
nearest centimeter so it would make sense to take $\delta=1$ cm and the
$\mu_{i}$ are consecutive integer values in centimeters and all values in
$[\mu_{i}-\delta/2,\mu_{i}+\delta/2)$ are considered effectively equivalent.
For the neutrino problem there is undoubtedly a guaranteed upper bound on the
mass. As discussed for Fieller's problem, priors are chosen via elicitation
and continuous priors are considered here with the previously described
discretization applied for computations when necessary. Results for two priors
are presented for comparison purposes.

The first prior $\pi_{1}$ is taken to be a beta$(\alpha_{0},\beta_{0})$
distribution on the interval $(l_{0},u_{0})$ with the elicitation procedure as
described in Evans, Guttman and Li (2017) although others are possible. For
this $\mu=l_{0}+(u_{0}-l_{0})z$ where $z\sim$ beta$(\alpha_{0},\beta_{0}).$
The values of $(\alpha_{0},\beta_{0})$ are specified as follows. First it is
required that $\alpha_{0},\beta_{0}\geq1$ to ensure unimodality and no
singularities. Next a proper subinterval $(l_{1},u_{1})\subset(l_{0},u_{0})$
is specified such that $\mu\in(l_{1},u_{1})$ with prior probability $\gamma.$
Typically $\gamma$ will be a large probability (like $0.99$ or higher)
reflecting the fact that $\mu\in(l_{1},u_{1})$ is known to be true with
virtual certainty. Then the mode $m_{0}$ is taken to be equal to a value in
$(l_{1},u_{1}),$ such as $m_{0}=(l_{1}+u_{1})/2,$ which implies $m_{0}%
=l_{0}+(u_{0}-l_{0})(\alpha_{0}-1)/\tau_{0}\in(l,u)$ where $\tau_{0}%
=\alpha_{0}+\beta_{0}-2.$ This leads to values for $\alpha_{0}\ $and
$\beta_{0}$ as%
\[
\alpha_{0}=\tau_{0}(m_{0}-l_{0})/(u_{0}-l_{0})+1,\beta_{0}=\tau_{0}%
(u_{0}-m_{0})/(u_{0}-l_{0})+1
\]
that are fully specified once $\tau_{0}$ is chosen. The value of $\tau_{0}$
controls the dispersion of the beta$(\alpha_{0},\beta_{0})$ and, with the cdf
denoted beta$(\cdot,\alpha_{0},\beta_{0}),$ we want beta$((u_{1}-l_{0}%
)/(u_{0}-l_{0}),\alpha_{0},\beta_{0})-$beta$((l_{1}-l_{0})/(u_{0}%
-l_{0}),\alpha_{0},\beta_{0})=\gamma$, and this is easily solved for $\tau
_{0}$ by an iterative procedure based on bisection. For example, with
$(l_{0},u_{0})=(0,10),(l_{1},u_{1})=(0.5,9.5),m_{0}=5$ and $\gamma=0.99$ this
leads to $(\alpha_{0},\beta_{0})=(2.20,2.20).$ Note that if $(l_{1}%
,u_{1})=(l_{0},u_{0}),$ then use the uniform prior, namely, $\alpha_{0}%
=\beta_{0}=1$ which is the noninformative case.

The second prior $\pi_{2}$ is taken to be a $N(\mu_{0},\tau_{0}^{2})$
constrained to the interval $(l_{0},u_{0}).$ The interval $(l_{1},u_{1})$ is
selected as before and $\mu_{0}\in(l_{1},u_{1})$ is \ specified while
$\tau_{0}$ is chosen to satisfy $(\Phi((u_{1}-\mu_{0})/\tau_{0})-\Phi
((l_{1}-\mu_{0})/\tau_{0}))/(\Phi((u_{0}-\mu_{0})/\tau_{0})-\Phi((l_{0}%
-\mu_{0})/\tau_{0}))=\gamma$ using the cdf of the prior. As $\tau
_{0}\rightarrow0$ the prior content of $(l_{1},u_{1})$ goes to 1 and as
$\tau_{0}\rightarrow\infty$ the limit, using L'H\^{o}pital, is $(u_{1}%
-l_{1})/(u_{0}-l_{0}).$ So provided $(u_{1}-l_{1})/(u_{0}-l_{0})\leq\gamma$
there is a solution for $\tau_{0}.$ For example, with $(l_{0},u_{0}%
)=(0,10),(l_{1},u_{1})=(0.5,9.5),m_{0}=5$ and $\gamma=0.99$ this leads to
$(\mu_{0},\tau_{0}^{2})=(5,1.92^{2}).$

The bias against hypothesis $H_{0}:\mu=\mu_{\ast}$ is given by $M(RB_{i}%
(\mu_{\ast}\,|\,\bar{x})\leq1\,|\,\mu_{\ast})$ where $M(\cdot\,|\,\mu_{\ast})$
is the $N(\mu_{\ast},\sigma_{0}^{2}/n)$ measure, $RB_{i}(\mu_{\ast}%
\,|\,\bar{x})=m_{i}(\bar{x}\,|\,\mu_{\ast})/m_{i}(\bar{x})$ and $m_{i}(\bar
{x})=(n/\sigma_{0}^{2})^{1/2}\int_{l_{0}}^{u_{0}}\varphi\left(  n^{1/2}%
(\bar{x}-\mu)/\sigma_{0}\right)  \pi_{i}(\mu)\,d\mu$ is the prior predictive
density of $\bar{x}$ using prior $\pi_{i}.$ The difficulty in evaluating
$M(RB_{i}(\mu_{\ast}\,|\,\bar{x})\leq1\,|\,\mu_{\ast})$ arises from the need
to evaluate $m_{i}(\bar{x})$ to obtain $RB_{i}(\mu_{\ast}\,|\,\bar{x})$ for
each $\bar{x}$\ generated from the $N(\mu_{\ast},\sigma_{0}^{2}/n)$
distribution$.$ For this we proceed via an approximation where a sample of $N$
is obtained from the $M_{i}$ distribution by generating $\mu\sim\pi_{i}%
,\bar{x}\,|\,\mu\sim N(\mu,\sigma_{0}^{2}/n),$ the interval $(l_{0}%
-3\sigma_{0}/\sqrt{n},u_{0}+3\sigma_{0}/\sqrt{n})$ is divided into $k$ equal
length subintervals and the proportion of $\bar{x}$ values falling in each of
the $k$ intervals is recorded. The probabilities $p_{1},\ldots,p_{k}$ of these
intervals with respect to the $N(\mu_{\ast},\sigma_{0}^{2}/n)$ are computed
and the relative belief ratios $RB(\mu_{\ast}\,|\,\bar{x})$ are then estimated
by the ratios of the $p_{i}$ to the relevant proportion obtained from sampling
from $M_{i}.$ Finally, the $p_{i}$ probabilities for the intervals where the
estimated relative belief ratio is less than or equal to $1$ are summed to
give the estimate of the bias against. Clearly as $k$ and $N$ increase this
approximation will converge to $M(RB_{i}(\mu_{\ast}\,|\,\bar{x})\leq
1\,|\,\mu_{\ast}).$ For these computations values of $N=10^{6}$ and of $k$ of
at least $10^{3},$ where the choice depended on $n,$ were used$.$ For example,
with $\sigma_{0}^{2}=1,$ the other constants as previously specified and
$\mu_{\ast}=4,$ Table \ref{t1} gives the values of the bias against for
different sample sizes and two different priors. It is seen that bias against
is not a problem with either prior.
\begin{table}[tbp] \centering
\begin{tabular}
[c]{|r|l|l|}\hline
$n$ & Bias against $H_{0}:\mu_{\ast}=4$ with $\pi_{1}$ & Bias against
$H_{0}:\mu_{\ast}=4$ with $\pi_{2}$\\\hline
$10$ & \multicolumn{1}{|c|}{$0.039$} & \multicolumn{1}{|c|}{$0.047$}\\
$20$ & \multicolumn{1}{|c|}{$0.026$} & \multicolumn{1}{|c|}{$0.032$}\\
$50$ & \multicolumn{1}{|c|}{$0.015$} & \multicolumn{1}{|c|}{$0.019$}\\
$100$ & \multicolumn{1}{|c|}{$0.011$} & \multicolumn{1}{|c|}{$0.013$}\\
$500$ & \multicolumn{1}{|c|}{$0.004$} & \multicolumn{1}{|c|}{$0.005$}\\\hline
\end{tabular}
\caption{Bias against values in Example 3  for testing $H_0:\mu_*=4$ for various sample sizes.}\label{t1}%
\end{table}%

\begin{figure}[t]%
\centering
\includegraphics[
height=2.2433in,
width=2.6368in
]%
{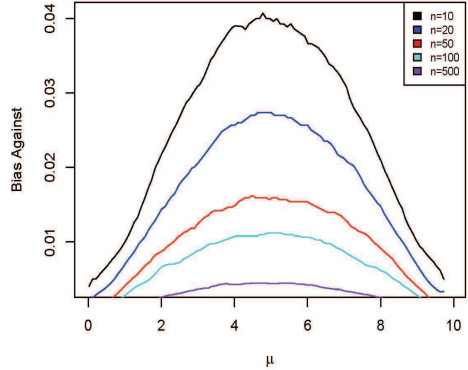}%
\caption{Bias against in Example 3 when using prior $\pi_{1}$ for various
sample sizes.}%
\label{mandelfig1}%
\end{figure}

\begin{figure}[t]%
\centering
\includegraphics[
height=2.1932in,
width=2.5797in
]%
{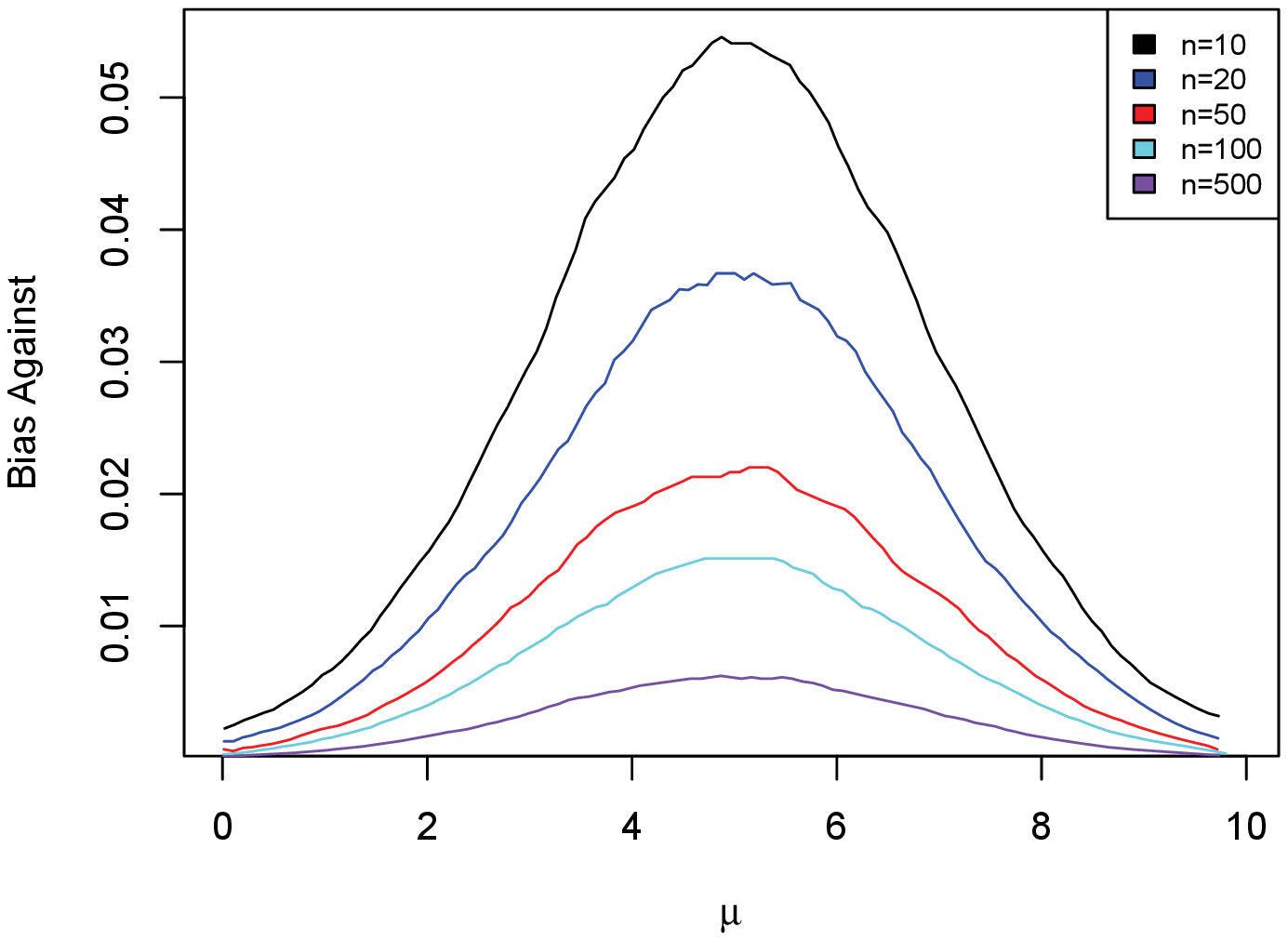}%
\caption{Bias against in Example 3 when using prior $\pi_{2}$ for various
sample sizes.}%
\label{mandelfig2}%
\end{figure}

Figure \ref{mandelfig1} is a graph of $M(RB_{i}(\mu\,|\,\bar{x})\leq
1\,|\,\mu)$ as a function of $\mu$ when using $\pi_{1}$ for various $n$ and
Figure \ref{mandelfig2} is the graph using $\pi_{2}.$ It is seen that the bias
against is maximized at a value $\mu_{\max}$ which implies that $M(RB_{i}%
(\mu_{\max}\,|\,\bar{x})\leq1\,|\,\mu_{\max})$ serves as an upper bound on the
bias against for estimation purposes. As such $M(RB_{i}(\mu_{\max}\,|\,\bar
{x})>1\,|\,\mu_{\max})$ is a lower bound on the coverage probabilities
$M(\mu\in Pl_{i}(\bar{x})\,|\,\mu)$ for $\mu\in(l_{0},u_{0})\ $where
$Pl_{i}(\bar{x})$ is the plausible region based on $\pi_{i}.$ As such
$Pl_{i}(\bar{x})$ as a confidence region has coverage probability
$M(RB_{i}(\mu_{\max}\,|\,\bar{x})>1\,|\,\mu_{\max})$ or greater. Table
\ref{t2} contains the confidence values for $Pl_{i}(\bar{x})$ for various
sample sizes. So it is seen that a $95\%$ frequentist coverage is achieved
fairly easily. It is to be noted, however, that the confidence is a priori and
the correct measure of belief that the true value is in $Pl_{i}(\bar{x}),$
based on the principle of conditional probability, is the posterior
probability. The Bayesian a priori coverage probability for $Pl_{i}(\bar{x})$
is also recorded in Table \ref{t2}. In many ways these coverage probabilities
can be considered as more appropriate than the pure frequentist coverage as
they take into account what is known about $\mu$ through the prior. There is
very little difference in this example.%
%

\begin{table}[tbp] \centering
\begin{tabular}
[c]{|r|l|l|}\hline
$n$ & $%
\begin{array}
[c]{c}%
\text{Confidence level of }Pl(\bar{x})\\
\text{using }\pi_{1}\text{ (Bayes)}%
\end{array}
$ & $%
\begin{array}
[c]{c}%
\text{Confidence level of }Pl(\bar{x})\\
\text{using }\pi_{2}\text{ (Bayes)}%
\end{array}
$\\\hline
$10$ & \multicolumn{1}{|c|}{$0.958$ $(0.969)$} & \multicolumn{1}{|c|}{$0.945$
$(0.961)$}\\
$20$ & \multicolumn{1}{|c|}{$0.973$ $(0.979)$} & \multicolumn{1}{|c|}{$0.962$
$(0.974)$}\\
$50$ & \multicolumn{1}{|c|}{$0.984$ $(0.988)$} & \multicolumn{1}{|c|}{$0.977$
$(0.985)$}\\
$100$ & \multicolumn{1}{|c|}{$0.988$ $(0.991)$} & \multicolumn{1}{|c|}{$0.985$
$(0.989)$}\\
$500$ & \multicolumn{1}{|c|}{$0.995$ $(0.997)$} & \multicolumn{1}{|c|}{$0.994$
$(0.996)$}\\\hline
\end{tabular}
\caption{Frequentist (Bayesian) confidence that $Pl_i(\bar{x})$ contains the true value in Example 3 for various sample sizes.}\label{t2}%
\end{table}%

\begin{figure}[ptbh]%
\centering
\includegraphics[
height=2.5374in,
width=2.9827in
]%
{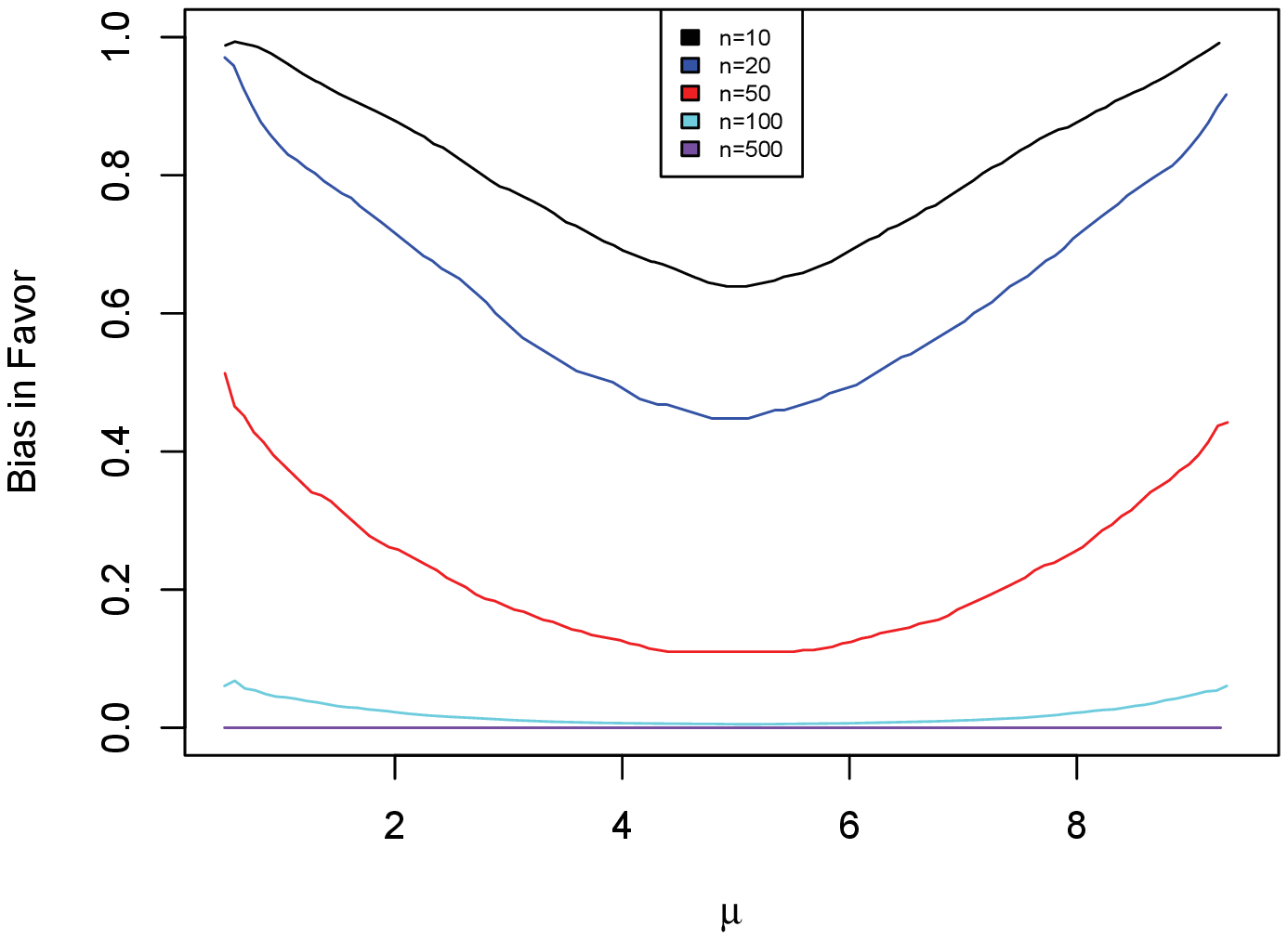}%
\caption{Bias in favor in Example 3 with prior $\pi_{2}$ for various $n$ and
$\delta.$}%
\label{mandelfig3}%
\end{figure}

It is also necessary to be concerned about bias in favor of $H_{0}:\mu
=\mu_{\ast}.$ Generally this bias is the more serious concern because of a
predilection towards the use of diffuse priors as these generally induce bias
in favor. Table \ref{t3} presents the bias in favor for different $n$ and
$\delta$ and Figure \ref{mandelfig3} is a plot of the bias in favor when using
$\pi_{2}$ and a similar plot is obtained with $\pi_{1}.$ Similar results are
obtained for the bias in favor for estimation as presented in Table \ref{t4}.
It is seen that the bias in favor in both problems can be substantial and for
a given prior and $\delta$ this can only be decreased by increasing the sample
size or by increasing $\delta$. One needs to be realistic about what accuracy
is necessary both, to make sure the study is returning results of sufficient
accuracy, and that resources are not being wasted. For estimation the bias in
favor is measured by averaging the bias in favor\ with respect to the prior.
As Table \ref{t4} indicates, large sample sizes are needed to make sure the
bias in favor is small, although this can be mitigated by taking $\delta$
larger.\smallskip\
%

\begin{table}[tbp] \centering
\begin{tabular}
[c]{|r|l|l|}\hline
$n$ & $%
\begin{array}
[c]{c}%
\text{Bias in favor of }H_{0}:\mu_{\ast}=4\\
\text{using }\pi_{1},\,\delta=0.5\,(0.1)
\end{array}
$ & $%
\begin{array}
[c]{c}%
\text{Bias in favor of }H_{0}:\mu_{\ast}=4\\
\text{using }\pi_{2},\,\delta=0.5\,(0.1)
\end{array}
$\\\hline
$10$ & \multicolumn{1}{|c|}{$0.700$ $(0.954)$} & \multicolumn{1}{|c|}{$0.686$
$(0.947)$}\\
$20$ & \multicolumn{1}{|c|}{$0.515$ $(0.962)$} & \multicolumn{1}{|c|}{$0.495$
$(0.957)$}\\
$50$ & \multicolumn{1}{|c|}{$0.136$ $(0.957)$} & \multicolumn{1}{|c|}{$0.123$
$(0.951)$}\\
$100$ & \multicolumn{1}{|c|}{$0.008$ $(0.941)$} & \multicolumn{1}{|c|}{$0.006$
$(0.934)$}\\
$500$ & \multicolumn{1}{|c|}{$0.000$ $(0.744)$} & \multicolumn{1}{|c|}{$0.000$
$(0.716)$}\\\hline
\end{tabular}
\caption{Bias in favor of $\mu_*$ in Example 3 for various sample sizes and meaningful differences.}\label{t3}%
\end{table}%
%

\begin{table}[tbp] \centering
\begin{tabular}
[c]{|r|l|l|}\hline
$n$ & $%
\begin{array}
[c]{c}%
\text{Bias in favor for estimation}\\
\text{using }\pi_{1},\,\delta=0.5\,(0.1)
\end{array}
$ & $%
\begin{array}
[c]{c}%
\text{Bias in favor for estimation}\\
\text{using }\pi_{2},\,\delta=0.5\,(0.1)
\end{array}
$\\\hline
$10$ & \multicolumn{1}{|c|}{$0.756$ $(0.965)$} & \multicolumn{1}{|c|}{$0.686$
$(0.947)$}\\
$20$ & \multicolumn{1}{|c|}{$0.569$ $(0.969)$} & \multicolumn{1}{|c|}{$0.495$
$(0.956)$}\\
$50$ & \multicolumn{1}{|c|}{$0.176$ $(0.966)$} & \multicolumn{1}{|c|}{$0.123$
$(0.951)$}\\
$100$ & \multicolumn{1}{|c|}{$0.012$ $(0.951)$} & \multicolumn{1}{|c|}{$0.006$
$(0.934)$}\\
$500$ & \multicolumn{1}{|c|}{$0.000$ $(0.763)$} & \multicolumn{1}{|c|}{$0.000$
$(0.716)$}\\\hline
\end{tabular}
\caption{Bias in favor for estimation in Example 3  for various sample sizes and meaningful differences.}\label{t4}%
\end{table}%

\noindent\textbf{Example 4}. \textit{Poisson with constrained mean.}

Suppose that count measurements $x_{1},\ldots,x_{n}$ are $i.i.d.$
Poisson$(\lambda)$ where $\lambda$ is known to lie in the interval
$(l_{0},u_{0}).$ The Poisson distribution arises as follows: suppose an event
occurs in a time interval of length 1 unit with probability $p$ and there are
$N$ independent opportunities for such events to occur. Then $\lambda\approx
Np$ for large $N$ and so $\lambda$ represents the rate at which the event
occurs in such a time interval. As discussed in Mandelkern (2002) sometimes
this rate is known to be at least $l_{0}>0$ and, as with the normal example,
without loss of generality, it will be supposed that $u_{0}$ is finite as
well. Again it is necessary to specify $\delta>0$ such that two values of
$\lambda$ that differ less than this are effectively equivalent and also
specify the grid of $\lambda_{i}$ values as was done as in Example 3.

Many possibilities exist for a prior $\pi$ but attention is restricted here to
a gamma$_{rate}(\alpha_{0},\beta_{0})$ prior. Again an interval $(l_{1}%
,u_{1})\subset(l_{0},u_{0})$ is specified together with a probability
$\gamma=\Pi((l_{1},u_{1}))$ and the mode $m_{0}\in(l_{1},u_{1}).$ Any value in
$(l_{1},u_{1})$ is allowed for the mode but the value $m_{0}=(l_{1}+u_{1})/2$
is selected here. Then $m_{0}=(\alpha_{0}-1)/\beta_{0}$ and $\pi$ is a
gamma$_{rate}(1+m_{0}\beta_{0},\beta_{0})$ with $\beta_{0}$ determined by
$\gamma$ which can be solved for iteratively using bisection. As a specific
example suppose $(l_{0},u_{0})=(3,10)$ and $(l_{1},u_{1})=(3.5,9.5).$ This
implies that the prior\ is a gamma$_{rate}(37.20,5.57)$ distribution.

The bias against $H_{0}:\lambda=6.2$ is recorded in Table \ref{t5} and for
modest sample sizes it is seen that this is well-controlled. Table \ref{t6}
provides the lower bounds on the confidence levels and the exact Bayesian
coverages for the plausible interval for various $n$ and $\delta.$ The
coverage probabilities are reasonable for $n\geq10.$%

\begin{table}[tbp] \centering
\begin{tabular}
[c]{|r|l|}\hline
$n$ & Bias against $H_{0}:\lambda=6.5$ with $\delta=0.5(1)$\\\hline
$1$ & \multicolumn{1}{|c|}{$0.287\,(0.286)$}\\
$10$ & \multicolumn{1}{|c|}{$0.193\,(0.180)$}\\
$20$ & \multicolumn{1}{|c|}{$0.085\,(0.045)$}\\
$50$ & \multicolumn{1}{|c|}{$0.045\,(0.011)$}\\
$100$ & \multicolumn{1}{|c|}{$0.001\,(0.000)$}\\
$500$ & \multicolumn{1}{|c|}{$0.000\,(0.000)$}\\\hline
\end{tabular}
\caption{Bias against values in Example 4 for testing  $H_0:\lambda =6.2$ for various sample sizes and meaningful differences.}\label{t5}%
\end{table}%
%

\begin{table}[tbp] \centering
\begin{tabular}
[c]{|r|l|l|}\hline
$n$ & $%
\begin{array}
[c]{c}%
\text{Confidence level of }Pl(\bar{x})\\
\text{using }\pi,\delta=0.5\text{ (Bayes)}%
\end{array}
$ & $%
\begin{array}
[c]{c}%
\text{Confidence level of }Pl(\bar{x})\\
\text{using }\pi,\delta=1.0\text{ (Bayes)}%
\end{array}
$\\\hline
$1$ & \multicolumn{1}{|c|}{$0.581\,(0.667)$} &
\multicolumn{1}{|c|}{$0.614\,(0.663)$}\\
$10$ & \multicolumn{1}{|c|}{$0.811$ $(0.840)$} & \multicolumn{1}{|c|}{$0.828$
$(0.858)$}\\
$20$ & \multicolumn{1}{|c|}{$0.843$ $(0.878)$} & \multicolumn{1}{|c|}{$0.865$
$(0.903)$}\\
$50$ & \multicolumn{1}{|c|}{$0.908$ $(0.935)$} & \multicolumn{1}{|c|}{$0.948$
$(0.966)$}\\
$100$ & \multicolumn{1}{|c|}{$0.950$ $(0.966)$} & \multicolumn{1}{|c|}{$0.986$
$(0.991)$}\\
$500$ & \multicolumn{1}{|c|}{$0.998$ $(0.999)$} & \multicolumn{1}{|c|}{$0.998$
$(1.000)$}\\\hline
\end{tabular}
\caption{Frequentist (Bayesian) confidence that $Pl_i(\bar{x})$ contains the true value in Example 4  for various sample sizes and meaningful differences.}\label{t6}%
\end{table}%

Table \ref{t7} presents values of the bias in favor of the hypothesis
$H_{0}:\lambda=6.2$ for various $n$ and $\delta.$ It is seen that there is
appreciable bias in favor of $H_{0}$ when $\delta=0.5$ unless $n\geq500.$ When
$\delta=1.0,$ however, much smaller sample sizes give reasonable values. Table
\ref{t8} provides values of for the bias in favor for estimation for various
$n$ and $\delta$ and large sample sizes are needed to get the bias in favor
down to acceptable levels.%

\begin{table}[tbp] \centering
\begin{tabular}
[c]{|r|l|}\hline
$n$ & $\text{Bias in favor of }H_{0}:\lambda=6.2$ $\text{using }\pi
,\delta=0.5\,(1.0)$\\\hline
$1$ & \multicolumn{1}{|c|}{$0.781\ (0.840)$}\\
$10$ & \multicolumn{1}{|c|}{$0.800$ $(0.673)$}\\
$20$ & \multicolumn{1}{|c|}{$0.667$ $(0.218)$}\\
$50$ & \multicolumn{1}{|c|}{$0.522\,(0.074)$}\\
$100$ & \multicolumn{1}{|c|}{$0.105$ $(0.000)$}\\
$500$ & \multicolumn{1}{|c|}{$0.016$ $(0.000)$}\\\hline
\end{tabular}
\caption{Bias in favor of $H_0 : \lambda=6.2$ in Example 4 for various sample sizes and meaningful differences.}\label{t7}%
\end{table}%
%

\begin{table}[tp] \centering
\begin{tabular}
[c]{|r|l|}\hline
$n$ & Bias in favor for estimation using $\pi,\delta=0.5\,(1.0)$\\\hline
$1$ & \multicolumn{1}{|c|}{$0.732\,(0.744)$}\\
$10$ & \multicolumn{1}{|c|}{$0.865$ $(0.778)$}\\
$20$ & \multicolumn{1}{|c|}{$0.850$ $(0.652)$}\\
$50$ & \multicolumn{1}{|c|}{$0.775$ $(0.404)$}\\
$100$ & \multicolumn{1}{|c|}{$0.662$ $(0.197)$}\\
$500$ & \multicolumn{1}{|c|}{$0.204$ $(0.002)$}\\\hline
\end{tabular}
\caption{Bias in favor for estimation in Example 4  for various sample sizes and meaningful differences.}\label{t8}%
\end{table}%

\section{Conclusions}

The approach taken here to the construction of regions, whether confidence or
credible, is somewhat different than what is typically done where a
probability $\gamma$ is stated, whether as a confidence or as a posterior
probability, and then the region is constructed based on the observed data and
this probability. Rather, using the principle of evidence, the plausible
region is obtained as consisting of those values for which there is evidence
in favor of them being the true value and then quoting the posterior
probability of the region as a measure of the degree of belief that the true
value is in the stated region. Confidence here is an a priori concept which
the experimenter\ uses, before the data is collected, to ensure the experiment
will lead to reliable results.

Mandelkern (2002) states five desiderata that an assessment of the accuracy of
an estimate via a confidence interval should satisfy. These are now stated
with an assessment of how well the methodology described here meets a
requirement.\smallskip

\noindent(i) \textit{Confidence bounds are determined using a well-defined
principle, which is neither arbitrary nor subjective. }The bounds stated here
are fully determined by the principle of evidence. This principle is universal
in the sense that it is applicable to all statistical problems and is not
tailored to problems with bounded parameters. No optimality criteria are
required although the regions obtained do have optimal properties.

\noindent(ii) \textit{They do not depend upon prior knowledge of the parameter
apart from its domain. }The principle of evidence requires that a proper prior
probability distribution be stated. It is to be noted, however, that the
methodology used here includes an elicitation algorithm for the choice of the
prior, the measurement and control of the bias induced by the prior-model
combination and the checking of the model and prior against the data to see if
they are contradicted. It is also the case that the check on the prior is a
check on any bounds assumed for the parameter. Objectivity is a necessary
aspect of scientific work although difficult to characterize precisely. For
example, frequentist methods are not objective as these involve subjective
choices made by a statistician. While objectivity is the ideal it is necessary
to recognize that, while it is unattainable, it can be approached via
methodologies that check subjective choices against the objective data and
measure and control the bias that such choices may induce.

\noindent(iii) \textit{They are equivariant under one-to-one transformation of
the data. }The inference methods described here are fully invariant under
\textit{all} smooth reparameterizations. The intervals are (integrated)
likelihood intervals but with the additional characteristic that these
intervals contain only values for which there is evidence in favor of being
the true value. Confidence, likelihood and credible intervals do not possess
this property.

\noindent(iv) \textit{They convey an estimate of the experimental uncertainty.
}In addition to their lengths, the intervals here satisfy a confidence
condition as well as providing the posterior probability that the true value
is in the interval. The confidence is seen as an a priori assessment of the
quality of the experiment while the posterior probability and the length of
the interval are assessments of the accuracy of the estimate based upon the
observed data. No matter what data is obtained, these intervals are never absurd.

\noindent(v) \textit{They correspond to a precise statement of probability.
}The Bayesian a priori coverage is precise and a precise (sharp) lower bound
is determined for the conditional, given the true value, prior coverage. These
coverages can be set a priori by choice of sample size. The intervals also
have a precise posterior probability which is the correct measure of belief
that the interval based on the observed data contains the true value. The
biases are a priori probabilities that provide an assessment of the quality of
an experiment. So, for example, if the a priori coverage is low, then this
suggests that the results have to be treated with caution even if the interval
is short and has a high posterior probability.

\section*{Acknowledgements}

John Edwards assisted with the computations concerning Mandelkern's examples
and was supported by a Natural Sciences and Engineering Research Council of
Canada Undergraduate Student Research Award for this work.

\section*{References}

\noindent Evans, M. (2015) Measuring Statistical Evidence Using Relative
Belief. Monographs on Statistics and Applied Probability 144, CRC Press,
Taylor \& Francis Group.\textbf{\smallskip}

\noindent Evans, M., Guttman, I. and Li, P. (2017) Prior elicitation,
assessment and inference with a Dirichlet prior. Entropy 2017, 19(10), 564;
doi:10.3390/e1910056.\textbf{\smallskip}

\noindent Evans, M. and Moshonov, H. (2006) Checking for prior-data conflict.
Bayesian Analysis, 1, 4, 893-914.\textbf{\smallskip}

\noindent Fieller, E. C. (1954) Some problems in interval estimation. JRSSB,
16, 2, 175-186. doi.org/10.1111/j.2517-6161.1954.tb00159.\textbf{\smallskip}

\noindent Fraser, D. A. S., Reid, N. and Lin, W. (2018) When should modes of
inference disagree? Some simple but challenging examples. Ann. Appl. Stat. 12
(2) 750 - 770, doi.org/10.1214/18-AOAS1160SF.\textbf{\smallskip}

\noindent Geary, R. C. (1930). The frequency distribution of the quotient of
two normal variates. J. Royal Statist. Soc. 97:442--446.
doi.org/10.2307/2342070\textbf{\smallskip}

\noindent Ghosh, M., Datta, G. S., Kim, D., and Sweeting, T. J. (2006)
Likelihood-based inference for the ratios of regression coefficients in linear
models. Ann. Inst. Stat. Math., 58: 457--473 DOI
10.1007/s10463-005-0027-3.\textbf{\smallskip}

\noindent Hinkley, D. V. (1969). On the ratio of two correlated normal random
variables. Biometrika 56:635-639.
doi.org/10.1093/biomet/57.3.683\textbf{\smallskip}

\noindent Mandelkern, M. (2002) Setting confidence intervals for bounded
parameters. Statistical Science 17(2): 149-172. doi:
10.1214/ss/1030550859\textbf{\smallskip}

\noindent Nott, D., Wang, X., Evans, M., and Englert, B-G. (2020) Checking for
prior-data conflict using prior to posterior divergences. Statistical Science,
35, 2, 234-253.\textbf{\smallskip}

\noindent Pham-Gia, T., Turkkan, N., and Marchand, E. (2006). Density of the
ratio of two normal random variables and applications. Communications in
Statistics. Theory and Methods, 35(9), 1569-1591.
doi.org/10.1080/03610920600683689\textbf{\smallskip}

\noindent Plante, A. (2020)\ A Gaussian alternative to using improper
confidence intervals. Canadian J. of Statistics, 48, 4, 773-801.
\end{document}